\documentclass[a4paper,11pt]{article}
\usepackage{amsmath,amsfonts,amssymb,amscd,amsthm,epsfig}
\usepackage[mathscr]{euscript}
\usepackage{hyperref}
\newtheorem{thm}{Theorem}[section]
\newtheorem{*thm}[thm]{Theorem}
\newtheorem{lemma}[thm]{Lemma}

\theoremstyle{definition}

\newtheorem{propr}[thm]{Property}
\newtheorem{proprs}[thm]{Properties}

\theoremstyle{remark}
 \newtheorem*{claim}{\textbf{Claim}}
\newtheorem*{rmq}{\textit{Remark}}
\newtheorem{rmk}[thm]{\textit{Remark}}
\renewcommand{\proof}{\noindent\textit{Proof}\/: \,\,}
\newcommand{\C}{{\mathbb{C}}}
\newcommand{\Q}{{\mathbb{Q}}}

\newcommand{\bV}{{\mathbb{V}}}

\newcommand{\bH}{{\mathbb{H}}}

\newcommand\CC{{\mathcal C}}

\newcommand{\II}{{\mathcal I}}
\newcommand{\comp}{\raise1pt\hbox{{$\scriptscriptstyle\circ$}}}
 
\newcommand{\underint}{\int_{\raise-4pt\hbox{\hskip-8pt $-$}}}
\newcommand{\overint}{\int^{\raise3pt\hbox{\hskip-7pt $-$}
}\hskip -4pt}

\newcommand\Tr{{}^{\mathsf{T}}\kern-0.9pt} 

\newcounter{lijstc}
{\end{list}}
{\end{list}}

\def\arrow(#1,#2)\dir(#3,#4)\long#5{\put(#1,#2){\vector(#3,#4){#5}}}

\newcommand\Grid{\setbox13=\vbox to 5\unitlength{\hrule width 109mm\vfill}
\setbox13=\vbox to 65mm{\offinterlineskip\leaders\copy13\vfill\kern-1pt\hrule}
\setbox14=\hbox to 5\unitlength{\vrule height 65mm\hfill}
\setbox14=\hbox to 109mm{\leaders\copy14\hfill\kern-2mm\vrule height 65mm}
\ht14=0pt\dp14=0pt\wd14=0pt \setbox13=\vbox to
0pt{\vss\box13\offinterlineskip\box14} \wd13=0pt\box13}

\def\mapright#1{\mathop{\vbox{\ialign{
                                ##\crcr
    ${\scriptstyle\hfil\;\;#1\;\;\hfil}$\crcr
 \noalign{\kern2pt\nointerlineskip}
    \rightarrowfill\crcr}}\;}}

\def\mapleft#1{\mathop{\vbox{\ialign{
                                ##\crcr
    ${\scriptstyle\hfil\;\;#1\;\;\hfil}$\crcr
 \noalign{\kern2pt\nointerlineskip}
    \leftarrowfill\crcr}}\;}}

\newcommand\rarrow[3]{\smash{\mathop{\hbox to#3{\rightarrowfill}}\limits
^{\scriptstyle#1}_{\scriptstyle#2}}}
 
\newcommand\larrow[3]{\smash{\mathop{\hbox to#3{\leftarrowfill}}\limits
^{\scriptstyle#1}_{\scriptstyle#2}}}

\renewcommand\setminus{-}
\def\into{\hookrightarrow}

\newcommand\im{\operatorname{Im}}

\newcommand\Hom{\mathop{\rm Hom}\nolimits}

\newenvironment{lijstb}{\vspace{-0.3ex}\begin{list}{\rm \Alph{lijstc})}{\usecounter{lijstc}%
\setlength{\leftmargin}{1ex}
\setlength{\rightmargin}{\leftmargin}}}%
{\end{list}}
\newcommand\mhm{\mathsf{MHM}}
\newcommand\perv{\mathop{\rm Perv}\nolimits}
\newcommand\rat{\mbox{\rm rat}}
\newcommand\Gr{\mathop{\rm Gr}\nolimits}
\newcommand\ic{\mathop{\II\CC}^\bullet\nolimits}

\title{Lowest Weights in Cohomology of Variations of Hodge Structure}
\author{Chris PETERS}
\begin{document}
\maketitle

\begin{abstract} Let $X$ be a smooth complex  projective variety,  let $j:U\into X$ an immersion of a Zariski open subset, and  let   $\bV$  be a variation of Hodge structure of weight $n$ over $U$. Then $IH^k(X, j_*\bV)$ is known to carry a pure Hodge structure of weigh $k+n$, while $H^k(U,\bV)$ carries a mixed Hodge structure of weight $\ge k+n$. In this note it is shown that the image of the natural map $IH^k(X,j_*\bV) \to H^k(U,\bV)$ is the lowest weight part of this mixed Hodge structure. The proof uses Saito's theory of mixed Hodge modules.
\end{abstract}
\section*{Introduction} 
For a smooth complex projective variety $X$ the  decomposition of complex valued differential $k$-forms into types  induces the Hodge decomposition for the De Rham   group  $H^k(X,\C)$  equipping this group with a   a pure weight $k$ Hodge structure.  For singular or non-compact  complex algebraic varieties this is no longer true in general. For instance $H^1(\C^*)$ has rank $1$ while it should have even rank if  it would carry a weight $1$ Hodge structure. Instead, cohomology groups  have  a \textbf{mixed Hodge structure}, i.e there is  a  rationally defined increasing weight filtration so that the $k$-th graded pieces carry a weight $k$  Hodge structure. In the example there is only one weight, namely $2$ and $H^1(\C^*)$ is pure of weight $(1,1)$.

Deligne \cite{Del71,Del74} constructed a good functorial theory for the cohomology of algebraic varieties and in the case of a smooth variety $U$ the weight filtration can be seen on the level of forms as follows.  First  choose a smooth compactification  $X$ such that $D=X\setminus U$ is a divisor with normal crossings. Cohomology of $U$ can now be calculated De Rham-style using rational forms having at most poles along $D$. The weight keeps track of the number of branches of $D$ along which one has actual poles. Lowest weight corresponds to forms that extend regularly across $D$. Indeed, $H^k(U)$ carries a mixed Hodge structure with  $W_{k-1}H^k(U)=0$ and where $W_kH^k(U)$ is the image of the restriction $H^k(X)\to H^k(U)$. 

If one studies how cohomology behaves under morphisms $f:Y \to X$ between algebraic varieties, the Leray spectral sequence tells one to look at $H^p(X,R^pf_*\underline{\Q}_Y)$. So  cohomology groups with varying coefficient system come up naturally.  Assume that $X$ and $Y$ are smooth projective. Then over the Zariski dense open subset $U\subset X$ of regular values of $f$  the sheaf $R^pf_*\underline{\Q}_Y$ is locally constant and   the fibers carry a  weight $p$ Hodge structure. In fact these can be assembled to give the prototype of what is called a \textbf{variation of weight $p$ Hodge structure} (cf. for instance \cite{CSP}). So it is then natural to look at $H^p(U,\bV)$ where $\bV$ is a local system which carries a weight $p$ variation of Hodge structure. It generalizes the previous case where $\bV=\underline{\Q}_U$ and as in that situation, it is known that there is a mixed Hodge structure on the cohomology group $H^p(U,\bV)$  provided   the local monodromy operators around  infinity are quasi-unipotent, a condition which is known  to hold in the geometric setting \cite{Sch73}.    The goal of this note is to show that  also in this setting, the lowest weight ``comes from the compactification''.

One should perhaps clarify what is meant by ``coming from the compactification'' because this is subtler than in the case of constant coefficients. Let me explain this for the case $\dim X=1$, following  \cite{Zuc79}. Let   $j:U\into X$ be  the embedding of $U$ into its compactification.  The sheaf $j_*\bV$  is quasi-isomorphic to the complex of holomorphic forms with values in $\bV$ and with $L^2$ growth conditions at the boundary (with respect to the Poincar\'e metric). Forgetting the growth conditions gives a complex which computes the cohomology of $\bV$ on $U$;  whence  a natural restriction map $H^p(X,j_*\bV) \to H^p(U, \bV)$. One of the main results from  \cite{Zuc79} states that the source has a pure weight $(p+q)$-Hodge structure which maps to the lowest weight part of a functorial mixed Hodge structure on the target.  As claimed before, in  the general situation this remains true  but it turns out to be easier  to replace  the approach using $L^2$-forms  by a purely topological approach, namely via  the intersection complex.  

In fact, for higher dimensional base one should work with $IH^p(X,j_*\bV)$. See Remark~\ref{L2andIC} where the two are compared.That  $IH^p(X,j_*\bV) $ carries a pure Hodge structure in the general situation (generalizing the curve case  \cite{Zuc79}) is  much harder. To tackle this problem M. Saito \cite{Sa88,Sa90} came up with  an elaborate  construction of mixed Hodge modules  in which  $D$-modules and perverse sheaves play an equally important role: in a mixed Hodge module they come glued in pairs related via the Riemann-Hilbert correspondence. It is crucial here  that the theory becomes functorial only after going over to the corresponding derived category. This functorial behavior makes it possible to use induction  reducing the general situation to the curve case where Zucker's result holds.

It is believed among experts that  the assertion about the lowest weights also holds, but in Saito's work this is not explicitly mentioned. I have written this note in order  to substantiate the expert's  believe. The result is stated as Theorem~\ref{mainresult}. The proof  turns out to be surprisingly subtle and uses in an essential way Saito's  functorial constructions, especially on the level of perverse sheaves. For that reason I start with a brief summary of the results I need  from this theory. Then I discuss the minimal ingredients from the theory of  mixed Hodge modules which are necessary in order to understand  the proof.  

I want to thank Stefan M\"uller-Stach for asking me this question and urging me to write down a proof.

\section{Perversity} \label{sec1} 
 I  won't recall here the rather technical definition  of perversity, but only give a minimal exposition to explain  the properties which will be used below.   I shall only be working with  the so called middle perversity which respects Poincar\'e duality. Full details can be found in \cite{B-B-D}.

Let $X$ be a complex algebraic variety of dimension $d$.    The category of  perverse ``sheaves''  of $\Q$-vector spaces on $X$, denoted by   $\perv (X;\Q)$, is an \emph{abelian category}.  The fact that it is abelian follows from its very construction as a core with respect to a $t$-structure. While the details of this are not so relevant for what follows, one needs to now that the starting point  is formed by  the \emph{constructible sheaves of $\Q$-vector spaces} on $X$.  By definition these are sheaves of $\Q$-vector spaces which  are locally constant on the strata  of some algebraic (and hence finite) stratification of $X$. The simplest   examples of such sheaves are the  locally constant sheaves on $X$ itself, or those that on some locally Zariski closed subset $Z$ of X are locally constant and zero elsewhere. The derived category of bounded complexes of such sheaves is also the bounded derived category of  $\perv (X;\Q)$. 

A core is defined  with respect to a so-called $t$-structure and in the perverse situation the $t$-structure is defined by the support and co-support conditions on complexes of constructible sheaves. By its very construction  a perverse sheaf is not a sheaf, but a complex of sheaves. 
\begin{rmk} \label{SupAndCosup} It is  important to have in mind is that the support condition for a perverse sheaf $F$ implies that $H^p(F)=0$ for $p>0$ while the co-support condition implies $H^p(F)=0$ for $p<-d$ (where $d=\dim X$):  perverse sheaves are complexes ``concentrated in degrees between $-d$ and $0$''.
\end{rmk}
On a smooth variety the support condition is even stronger and here a perverse sheaf is entirely  concentrated in degree $-d$. For instance the constant sheaf  on a smooth variety $X$ can be made perverse by placing it in degree $-d$: the complex $\underline{\Q}_X[d]$ is a perverse sheaf.  If $X$ is no longer smooth this complex has to be replaced by the so-called intersection complex $\ic_X$. More generally, if $U\subset X$ is a dense open subset of $X$ which consists of smooth points and $\bV$ is a (finite rank) local system of $\Q$-vector spaces on $U$ one can form  the  \textbf{intersection complex} $\ic_X(\bV[d])$\footnote{Some people write   $\ic_X(\bV)$ instead of $\ic_X(\bV[d])$} ; it  is the unique perverse extension of $\bV[d]$ to $X$.   By definition, its hypercohomology groups are the intersection cohomology groups:
\begin{equation}\label{eqn:thisisintcoh}
IH^k(X,\bV[d]):= \bH^k(\ic_X(\bV[d])).
\end{equation}

\begin{rmq} Even if $X$ itself is smooth  an intersection complex on $X$ need not be of the form  $\widetilde\bV[d]$ for some local system $\widetilde\bV$ on $X$ because of non-trivial monodromy ``around  infinity'' $X\setminus U$.\end{rmq}
The following two results  explain  the role  of these intersection complexes.
\begin{thm}[\protect{\cite[Chap.V, 4]{Bor}}]  \label{ICandPerv} Let $X$ be a  $d$-dimensional complex algebraic variety and let $U$ be a dense smooth subset of $X$ on which there is  a  a local system $\bV$ of finite dimensional $\Q$ vector spaces.   The intersection complex $\ic_X(\bV[d])$ is up to quasi-isomorphism the unique     complex of sheaves of $\Q$-vector spaces on $X$  which is perverse on $X$, which restricts over $U$ to $\bV[d]$ and which has no non-trivial perverse sub or quotient objects supported on $X\setminus U$.  
\end{thm}  
\begin{thm}[\cite{B-B-D}]  \label{SimplePerverseComplex} If $X$ is complex algebraic, $\perv(X;\Q)$ is artinian.   Its simple objects    are the intersection complexes  $F=\ic_Z(\bV[\dim Z])$ supported on an irreducible subspace $Z\subset X$  and  where $\bV$ is associated to an irreducible representation of $\pi_1(U)$, $U\subset Z$  the largest open subset of $Z$ over which $F$ is locally
constant.  
\end{thm}

\section{Mixed Hodge Modules} \label{sec2}
 
In this section I put together some properties of mixed Hodge modules which will be used in the sequel.
These properties are proven in \cite{Sa88} and \cite{Sa90}. I follow  largely the exposition \cite[Cha. 14]{PS} where mixed Hodge Modules are introduced axiomatically.  \medskip

Let $X$ be an algebraic variety. There exists an abelian category 
${\mhm}(X)$, the category of \textbf{mixed Hodge modules} on $X$  with the following
\begin{proprs}  \label{Props}
\begin{lijstb}
\item There is a faithful functor
\begin{equation}
\rat_X:D^{\rm b}\mhm(X)\to D^{\rm b}\perv(X;\Q).
\label{eqn:Rat}
\end{equation}
such that  $\mhm(X)$ corresponds to $\perv(X;\Q)$.  One  says  that $\rat_X M$ is the underlying
rational perverse sheaf of $M$. Moreover, let me  say that
\[
M\in \mhm(X) \mbox{ \textbf{is supported on }}   Z \iff  \rat_X M \, \mbox{is
supported on } Z.
\]

\item The category of mixed Hodge modules supported on a point is the category
of graded polarizable rational mixed Hodge structures; the functor ``$\rat$'' associates to the mixed
Hodge structure the underlying rational   vector space.
\item   Each object $M$ in ${\mhm}(X)$ admits a \textbf{weight filtration} $W$ such that
\begin{itemize}
\item morphisms preserve the weight filtration strictly;
\item the object $\Gr^W_kM$ is semisimple in $\mhm(X)$;
\item if $X$ is a point the $W$-filtration is the usual weight filtration for the mixed
Hodge structure.
\end{itemize}
Since $\mhm(X)$ is an abelian category, the cohomology groups of any complex of mixed
Hodge modules on $X$ are again  mixed Hodge modules on $X$. With this in mind, we say that
for complex $M^\bullet \in D^{\rm b}\mhm(X)$ the \textbf{weight} satisfies
\[ \mbox{weight}[M^\bullet] 
\left\{ \begin{array}{ll} 
\le n, \\
\ge n
\end{array} \right. \iff  \Gr_i^WH^j( M^\bullet)=0 \,\,
\left\{ \begin{array}{ll} 
\mbox{for } i>j+n\\
\mbox{for } i<j+n.  
\end{array} \right.  
\]
We observe  that  
 if we consider the  weight filtration on the mixed Hodge modules which constitute a complex  $M^\bullet\in D^{\rm b}\mhm(X)$ of mixed Hodge modules we   get a filtered complex in this category. 
\item For each morphism $f:X\to Y$ between complex algebraic varieties, there  are
  induced functors $f_* :D^{\rm b}\mhm(X)\to D^{\rm b}\mhm(Y)$ and $f^* :D^{\rm
b}\mhm(Y)\to D^{\rm b}\mhm(X)$   which lift the functors $Rf_*$ and $f^{-1}$  existing on the level of constructible complexes.

\item The functor  $f^*$ does   not increase weights in the sense that if $M^\bullet$
has weights $\le n$, the same is true for  $f^*M^\bullet$.
\item The functor  $f_*$ does  not decrease weights in the sense that if $M^\bullet$ has
weights $\ge n$,  the same is true for  $f_*M^\bullet$.

\item If $f$ is proper, $f_*$ preserves weights, i.e. $f_*$ neither increases nor decreases weights.
\end{lijstb}
\end{proprs}
These properties imply already various basic properties of mixed Hodge modules:
If $M^\bullet$ is a complex  of mixed Hodge modules on $X$ its cohomology $H^q M^\bullet $ is a mixed Hodge module on $X$. A consequence of Property  A) then is:

\begin{lemma}\label{WorkWithPerv} The cohomology functors $H^q:D^{\rm b}\mhm(X) \to \mhm(X)$
are compatible with the functor $\rat_X$ in the sense that for any bounded complex
$M^\bullet$ of mixed Hodge modules we have
\[
\rat_X [H^q M^\bullet]= \, ^\pi \! H^q[\rat_X M^\bullet],
\]
where the \emph{perverse} cohomology functor is used (see \S~\ref{sec1}).
\end{lemma}
Properties E) and B) imply:
\begin{lemma}  \label{HyperCohForMHS} Let  $a_X:X\to \mbox{\rm pt}$ be the constant
map  to the  point. For any complex $M^\bullet$ of mixed Hodge modules on $X$ 
\begin{equation} \label{eqn:HyperCohForMHS}
\bH^p(X,M^\bullet) := H^p((a_X)_* M^\bullet)
\end{equation}
is a  mixed Hodge structure. 
\end{lemma}

\section{Polarizable Variations of Hodge Structure} \label{sec3}

In this section $X$  is  a \emph{smooth} complex projective variety of dimension $d$ and   $j:U\into X$ is the inclusion of a smooth Zariski-dense open subset such the following holds:
\newline Assumption 1) \textit{$X\setminus U=D$ is a divisor with strict normal crossings.}
\newline
Recall  \cite[Th.~5.4.3]{Sa88}:
\begin{thm} \label{PolHS} Suppose that  $\bV$ is a polarizable variation on $U$ of weight $n$. If assumption 1) holds,  there is a pure weight $n$ mixed Hodge module $V^{\rm Hdg}$ on $U$  whose underlying perverse component is $\bV[d]$. 
\end{thm}
Pure weight $n$ mixed Hodge modules  are  exactly the \textbf{polarizable weight $n$ Hodge modules}. These form a semi-simple category (by property~\ref{Props}.C) and they satisfy moreover  the strict support condition: 
\begin{propr} A polarizable weight $n$ Hodge module  is a direct sum of polarizable weight $n$ Hodge modules which have strict support in some irreducible subvariety of $X$. \footnote{$M$  is said to have  \textbf{strict support in $Z$} if it is supported on $Z$  but  no quotient  or sub object  of $M$ has support  on a proper subvariety of $Z$.}
\end{propr}
Let me  make a second assumption on $\bV$:
\newline Assumption 2) \textit{The local monodromy operators around $D$ are quasi-unipotent.}
\newline
Then, by   \cite[3.20, 3.21]{Sa90} one has:
\begin{thm} \label{PolHM} If assumption 1) and 2) hold,  and if $\bV$ underlies a polarized  variation of Hodge  structures of weight $n$ on $U$, then there is a unique  
polarizable Hodge module $V^{\rm Hdg}_X$ of weight $n+d$ on $X$   having strict support in $X$ and which  restricts over $U$  to  $V^{\rm Hdg}$. 
\end{thm}
\begin{rmq} Note that this checks with the assertion in Theorem~\ref{ICandPerv} which holds for the rational component of the mixed Hodge modules.
\end{rmq}
The hypercohomology groups $\bH^k(X,\ic_X(\bV[d]))$ carry pure Hodge structures of weight $k+d+n$. This follows from the properties mentioned in \ref{sec2}:  consider the proper map $a_X:X \to $ pt. Then $[a_X]_*V^{\rm Hdg}$ is a \emph{complex} of pure mixed Hodge modules of weight $n+d$  over   a point  and its $k$-th cohomology has weight $k+d+n$. Since by \eqref{eqn:thisisintcoh} one has $IH^k(X,\bV)= \bH^k(X,\ic_X\bV)$, it follows that $IH^k(X,\bV)$ carries a pure  Hodge structure of weight $k+n$. 
\begin{rmk} \label{L2andIC}
By  \cite[Theorem 1.5]{CKS}    in this case the latter cohomology group can be identified with $L^2H^k(X,\bV)$ provided one measures integrability with respect to the  Poincar\'e metric  around infinity (one is in the normal crossing situation, so locally around infinity you have a product of disks and punctured disks).
Summarizing:
\[
\bH^k(\ic_X(\bV)) = IH^k(X, \bV)= L^2H^{k}(X ,\bV)  
\]
has a pure Hodge structure of weight  $k+n$.
\end{rmk}

Next I want to relate intersection cohomology and ordinary cohomology.  Start with the  very general situtation of a morphism  $f:X\to Y$  and a mixed Hodge module  $M$ on $Y$. The adjunction morphism $f^\#: M  \to f_* f^* M$ is a  morphism of complexes of mixed Hodge modules.  For any bounded complex $K^\bullet$ of sheaves on $X$, the identity $a_X=a_Y\comp f$ induces  a canonical identification $\bH^n(Y,  f_* K^\bullet)= \bH^n(X, K^\bullet)$. In particular this holds for  $K=f^*M$. Adjunction thus induces   a morphism  of mixed Hodge structures $ H^kf^\#: \bH^k( Y, M)  \to \bH^k(X, f^* M)$.
Specializing  the above to the inclusion $j:U\into X$    and  to $M:=V_X^{\rm Hdg}$, by the above remarks, the adjunction
\[
j^\#: V_X^{\rm Hdg} \to V^{\rm hdg}\]
induces    a morphism of mixed Hodge structures
 \begin{equation}\label{eqn:WhereItComesFrom}
H^kj^\#:  IH^k(X, \bV)   \to  H^k( U, \bV).
 \end{equation}

\begin{thm}  \label{mainresult} The image of $H^kj^\#$ under \eqref{eqn:WhereItComesFrom} is exactly the lowest weight part of $H^k( U, \bV)$
\end{thm}
\proof
By Property~\ref{Props}.G  the  complex $j_* V^{\rm Hdg}$ has weight $\ge n+d$.  The proposition follows from 
\begin{claim} The morphism $j^\#$ is injective and identifies  $V_X^{\rm Hdg}$ with the lowest weight part of $j_* V^{\rm Hdg}$.
\end{claim}
Suppose that this claim has been proven. Note that the map $H^kj^\#$ is the composition of the two maps $H^k(\mbox{\rm pt},  (a_X)_* V_X^{\rm Hdg}) \to  H^k(\mbox{\rm pt},  (a_X)_* j^\#V_X^{\rm Hdg})$ and $\iota_k: H^k(\mbox{\rm pt},  (a_X)_* j^\#V_X^{\rm Hdg}) \to H^k(\mbox{\rm pt},  (a_X)_* j_*V^{\rm Hdg})$. The first map is an isomorphism by the first part of the claim. By the second part of the claim, $j^\#V_X^{\rm Hdg}=W_{n+d}   j_*V^{\rm Hdg}$.
A priori $ j_*V^{\rm Hdg}$ is a complex of Hodge modules but in fact,  in this special case, since  the complement of $D$ is  a normal crossing divisor, Saito's construction of $j_*$ guarantees  that one obtains a single mixed Hodge module  $M$ (see \cite[94.2.11)]{Sa90}.  Then, 
by definition   $W_{n+d}H^{k}(\mbox{\rm pt},  (a_X)_* M)$ is  the image of the natural map
\[
  H^k(\mbox{\rm pt},  (a_X)_* W_{n+d} M)  \to H^k(\mbox{\rm pt},  (a_X)_*M)
\]
which is  the image of  the second map $\iota_k$. It follows that $W_{n+d}H^{k}  (a_X)_* M)$  is the image of $H^kj^\#$  which completes the proof.

 It remains to prove the claim.  I    use Theorem~\ref{PolHM} stating that the intersection complex for $\bV[d]$ underlies a mixed Hodge module which is  the unique polarisable Hodge module on $X$ having the property that it restricts to $V^{\rm Hdg}$ on $U$ and has no quotient  nor sub object supported on $D=X\setminus U$. It suffices therefore to show that the  lowest weight part $W_{d+n}M$ of $M=j_*V^{\rm Hdg}$ has these properties. It is  a  pure weight mixed Hodge module, and hence, by Property~\ref{Props}.C  a semi-simple object in the category of mixed Hodge modules. By construction, it restricts to $V^{\rm Hodge}$ on $U$.  By semi-simplicity a quotient object is also a sub object and hence it suffices to show that there are  no sub  objects $N$ supported on $D$, i.e. I need to show that $\Hom (N,  W_{n+d} M)=0$ in the \emph{derived category} of bounded complexes of mixed Hodge modules. This follows as soon as one shows the vanishing of   $\Hom (N,  M)$. Let $i:D\into X$ be the defining embedding. 
 Since $N=i_* i^*N$ by adjunction one has  $\Hom (N, M)=\Hom(i^*N, i^*M) $ and it suffices to prove $\Hom(i^*N, i^*M) =0$.
  
 Since by Property~\ref{Props}.A the perverse part is represented faithfully in the Hodge module  it suffices to show vanishing  of the latter group on the level of the rational components. 
The rational component of   $ M=j_*V^{\rm Hdg}$ is $j_*\bV[d]$. The constructible  sheaf $\widetilde \bV:=j_*\bV$ is   supported  as a locally constant sheaf on $U\cup D'$, where $D'$ is the union of components of $D$ along which the monodromy of $\bV$ is the identity.  Clearly $i^{-1} \widetilde\bV$ is supported as a locally constant sheaf  on $D'\subset D$. The rational component of $i^*j_*V^{\rm Hdg}$ is equal to $i^{-1} \widetilde\bV[d]$ and hence totally concentrated in degree $-d$. On the other hand the rational component of $i^*N$ is a perverse sheaf on $D$. By Remark~\ref{SupAndCosup}, since $\dim D=d-1$, this is  a complex $N^\bullet$ of constructible sheaves on $D$  with cohomology in the range $[-d+1,\dots,0]$. 

Now I recall some  elementary facts  from cohomological algebra.  Given a complex $K^\bullet$  in some abelian category  its two \textbf{truncations} are
given by
\begin{eqnarray*}
\tau_{\le k}K^\bullet & := &  \cdots K^{k-2}\to K^{k-1}\to \ker d^k\to 0\to
\cdots \\
\tau_{\ge k}K^\bullet  & := & \cdots \to 0\to \im d^{k-1}\to K^{k+1}\to
K^{k+2}\to \cdots \ . 
\end{eqnarray*}
These fit into a short exact sequence
\begin{equation}\label{eqn:Truncations}
0 \to \tau_{\le k} K^\bullet \to K^\bullet \to \tau_{\ge k+1}(K^\bullet)\to 0.
\end{equation}
Using \eqref{eqn:Truncations} one easily sees that $K^\bullet$ has the same cohomology in  degrees $\le k$   as its upper truncation $\tau_{\le k}K^\bullet$   while it has the same cohomology in degrees $\ge k$ as its lower truncation  $\tau_{\ge k} K^\bullet$.

In the case at hand, $N^\bullet$ has  no cohomology in the ranges $k>0$ and $k<-d+1$  and hence is quasi-isomorphic to the doubly  truncated complex $\tau_{\ge -d} \tau_{\le 0}  N^\bullet $ which is concentrated in degrees $[-d+1,0]$. In the derived category one thus has    $\Hom( N^\bullet, i^{-1}\widetilde{\bV} [d])=  \Hom(\tau_{\ge -d} \tau_{\le 0}  N^\bullet  , i^{-1} \widetilde\bV[d])=0$, as required.
\qed \endproof


\begin{thebibliography}{SaitoZuBB}
 
 \bibitem[B-B-D]{B-B-D} Beilinson, A., J. Bernstein and P.~Deligne:
Faisceaux pervers, in: \textsl{Analyse et topologie sur les espaces singuliers
I}, {Ast{\'e}risque}, \textbf{100}, (1982) 

 
 \bibitem[Bor84]{Bor} Borel, A. et  al.: \textsl{Intersection
cohomology}, Progress in Math. \textbf{50}, Birkh{\"a}user Verlag (1984)


 \bibitem[CKS]{CKS}Cattani, E. , A. Kaplan, and W. Schmid: $L^2$ and intersection cohomologies for a polarizable variation of Hodge structures. Inv. Math. \textbf{87}, 217--252 (1987).
 
 \bibitem[CSP]{CSP} Carlson, J., S. M\"uller-Stach, C. Peters: \textsl{Period mappings and Period Domains}, Cambr. stud. in Adv. Math. \textbf{85} Cambr. Univ. Press 2003.

\bibitem[Del71]{Del71} Deligne, P.: Th{\'e}orie de
Hodge II, Publ. Math. I.H.E.S, \textbf{40}, 5--58 (1971)

\bibitem[Del74]{Del74} Deligne, P.: Th{\'e}orie de Hodge III, Publ.
Math., I.~H.~E.~S, \textbf{44}, 5-77 (1974)

 \bibitem[PS]{PS} Peters, C. and J. Steenbrink: \textsl{Mixed Hodge Structures}, to appear in Ergebn. der Math. Wiss. Springer Verlag 2008.
 
 \bibitem[Sa88]{Sa88} Saito, M.: Modules de Hodge polarisables, Publ. RIMS. Kyoto
Univ. {\bf 24} 849--995 (1988)

\bibitem[Sa90]{Sa90} Saito, M.: Mixed Hodge Modules, Publ. Res. Inst. Math.
Sci.  \textbf{26}  221--333 (1990) 

\bibitem[Schm]{Sch73} Schmid, W.: Variation of Hodge structure: the
singularities of the period mapping, Invent. Math., \textbf{22}, 211--319 (1973) 


\bibitem[Zuc]{Zuc79} Zucker, S.: Hodge theory with degenerating
coefficients: $L_2$-cohomology in the Poincar\'e-metric, Ann. Math.,
\textbf{109}, 415--476 (1979)

 \end{thebibliography}
\end{document}